\documentclass[11pt,english]{article}
\usepackage[T1]{fontenc}
\usepackage[latin9]{inputenc}
\usepackage{geometry}
\usepackage{prettyref}
\usepackage{amsmath}
\usepackage{amsthm}
\usepackage{amssymb}
\usepackage{setspace}
\usepackage{esint}
\makeatletter

\newrefformat{chap}{Chapter \ref{#1}}
\newrefformat{prop}{Proposition \ref{#1}}
\newrefformat{cor}{Corollary \ref{#1}}
\newrefformat{rem}{Remark \ref{#1}}
\newrefformat{ex}{Example \ref{#1}}
\newrefformat{fig}{Figure \ref{#1}}

\allowdisplaybreaks

\makeatother

\usepackage{babel}
\begin{document}

\title{New Metric and Connections\\
in Statistical Manifolds\thanks{This work was partially funded by CNPq (Proc. 309055/2014-8).}}

\author{Rui F.\ Vigelis\thanks{Computer Engineering, Campus Sobral, Federal University of Ceará,
Sobral-CE, Brazil,\break rfvigelis@ufc.br.},\and David C.\ de Souza\thanks{Federal Institute of Ceará, Fortaleza-CE, Brazil, davidcs@ifce.edu.br.},
\and Charles C.\ Cavalcante\thanks{Wireless Telecommunication Research Group, Department of Teleinformatics
Engineering, Federal University of Ceará, Fortaleza-CE, Brazil, charles@ufc.br.}}
\maketitle
\begin{abstract}
We define a metric and a family of $\alpha$-connections in statistical
manifolds, based on $\varphi$-divergence, which emerges in the framework
of $\varphi$-families of probability distributions. This metric and
$\alpha$-connections generalize the Fisher information metric and
Amari's $\alpha$-connections. We also investigate the parallel transport
associated with the $\alpha$-connection for $\alpha=1$.
\end{abstract}

\section{Introduction}

In the framework of $\varphi$-families of probability distributions
\cite{Vigelis:2013}, the authors introduced a divergence $\mathcal{D}_{\varphi}(\cdot\mathbin{\Vert}\cdot)$
between probabilities distributions, called \textit{$\varphi$-divergence},
that generalizes the Kullback--Leibler divergence. Based on $\mathcal{D}_{\varphi}(\cdot\mathbin{\Vert}\cdot)$
we can define a new metric and connections in statistical manifolds.
The definition of metrics or connections in statistical manifolds
is a common subject in the literature \cite{Amari:2011,Amari:2012,Matsuzoe:2014}.
In our approach, the metric and $\alpha$-connections are intrinsically
related to $\varphi$-families. Moreover, they can be recognized as
a generalization of the Fisher information metric and Amari's $\alpha$-connections
\cite{Amari:2000,Calin:2014}.

Statistical manifolds are equipped with the Fisher information metric,
which is given in terms of the derivative of $l(t;\theta)=\log p(t;\theta)$.
Another metric can be defined if the logarithm $\log(\cdot)$ is replaced
by the inverse of a $\varphi$-function $\varphi(\cdot)$ \cite{Vigelis:2013}.
Instead of $l(t;\theta)=\log p(t;\theta)$, we can consider $f(t;\theta)=\varphi^{-1}(p(t;\theta))$.
The manifold equipped with this metric, which coincides with the metric
derived from $\mathcal{D}_{\varphi}(\cdot\mathbin{\Vert}\cdot)$,
is called a \emph{generalized statistical manifold}. 

Using the $\varphi$-divergence $\mathcal{D}_{\varphi}(\cdot\mathbin{\Vert}\cdot)$,
we can define a pair of mutually dual connections $D^{(1)}$ and $D^{(-1)}$,
and then a family of $\alpha$-connections $D^{(\alpha)}$. The connections
$D^{(1)}$ and $D^{(-1)}$ corresponds to the exponential and mixture
connections in classical information geometry. For example, in parametric
$\varphi$-families, whose definition is found in \prettyref{sec:phi-family},
the connection $D^{(1)}$ is flat (i.e, its torsion tensor $T$ and
curvature tensor $R$ vanish identically). As a consequence, a parametric
$\varphi$-family admits a parametrization in which the Christoffel
symbols $\Gamma_{ijk}^{(-1)}$ associated with $D^{(-1)}$ vanish
identically. In addition, parametric $\varphi$-families are examples
of Hessian manifolds \cite{Shima:2007}.

The rest of the paper is organized as follows. In \prettyref{sec:gen_stat_manifold},
we define the generalized statistical manifolds. \prettyref{sec:phi-family}
deals with parametric $\varphi$-families of probability distribution.
In \prettyref{sec:alpha_connections}, $\alpha$-connections are introduced.
The parallel transport associated with $D^{(1)}$ is investigated
in \prettyref{sec:parallel_transport}.

\section{Generalized Statistical Manifolds\label{sec:gen_stat_manifold}}

In this section, we provide a definition of generalized statistical
manifolds. We begin with the definition of $\varphi$-functions. Let
$(T,\Sigma,\mu)$ be a measure space. In the case $T=\mathbb{R}$
(or $T$ is a discrete set), the measure $\mu$ is considered to be
the Lebesgue measure (or the counting measure). A function $\varphi\colon\mathbb{R}\rightarrow(0,\infty)$
is said to be a \textit{$\varphi$-function} if the following conditions
are satisfied:
\begin{itemize}
\item [(a1)] $\varphi(\cdot)$ is convex,
\item [(a2)] $\lim_{u\rightarrow-\infty}\varphi(u)=0$ and $\lim_{u\rightarrow\infty}\varphi(u)=\infty$.
\end{itemize}
Moreover, we assume that a measurable function $u_{0}\colon T\rightarrow(0,\infty)$
can be found such that, for each measurable function $c\colon T\rightarrow\mathbb{R}$
such that $\varphi(c(t))>0$ and $\int_{T}\varphi(c(t))d\mu=1$, we
have
\begin{itemize}
\item [(a3)] ${\displaystyle \int_{T}\varphi(c(t)+\lambda u_{0}(t))d\mu<\infty}$,
for all $\lambda>0$.
\end{itemize}

The exponential function and the Kaniadakis' $\kappa$-exponential
function \cite{Kaniadakis:2002} satisfy conditions (a1)--(a3) \cite{Vigelis:2013}.
For $q\neq1$, the $q$-exponential function $\exp_{q}(\cdot)$ \cite{Tsallis:1994}
is not a $\varphi$-function, since its image is $[0,\infty)$. Notice
that if the set $T$ is finite, condition (a3) is always satisfied.
Condition (a3) is indispensable in the definition of non-parametric
families of probability distributions \cite{Vigelis:2013}.

A generalized statistical manifold is a family of probability distributions
$\mathcal{P}=\{p(t;\theta):\theta\in\Theta\}$, which is defined to
be contained in 
\[
\mathcal{P}_{\mu}=\biggl\{ p\in L^{0}:p>0\text{ and }\int_{T}pd\mu=1\biggr\},
\]
where $L^{0}$ denotes the set of all real-valued, measurable functions
on $T$, with equality $\mu$-a.e. Each $p_{\theta}(t):=p(t;\theta)$
is given in terms of parameters $\theta=(\theta^{1},\dots,\theta^{n})\in\Theta\subseteq\mathbb{R}^{n}$
by a one-to-one mapping. The family $\mathcal{P}$ is called a \emph{generalized
statistical manifold} if the following conditions are satisfied:
\begin{itemize}
\item [(P1)] $\Theta$ is a domain (an open and connected set) in $\mathbb{R}^{n}$. 
\item [(P2)] $p(t;\theta)$ is a differentiable function with respect to
$\theta$. 
\item [(P3)] The operations of integration with respect to $\mu$ and differentiation
with respect to $\theta^{i}$ commute.
\item [(P4)] The matrix $g=(g_{ij})$, which is defined by 
\begin{equation}
g_{ij}=-E_{\theta}'\Bigl[\frac{\partial^{2}f_{\theta}}{\partial\theta^{i}\partial\theta^{j}}\Bigr],\label{eq:g_ij_gen}
\end{equation}
is positive definite at each $\theta\in\Theta$, where $f_{\theta}(t)=f(t;\theta)=\varphi^{-1}(p(t;\theta))$
and 
\[
E_{\theta}'[\cdot]=\frac{\int_{T}(\cdot)\varphi'(f_{\theta})d\mu}{\int_{T}u_{0}\varphi'(f_{\theta})d\mu}.
\]

\end{itemize}

Notice that expression \prettyref{eq:g_ij_gen} reduces to the Fisher
information matrix in the case that $\varphi$ coincides with the
exponential function and $u_{0}=1$. Moreover, the right-hand side
of \prettyref{eq:g_ij_gen} is invariant under reparametrization.
The matrix $(g_{ij})$ can also be expressed as
\begin{equation}
g_{ij}=E_{\theta}''\Bigl[\frac{\partial f_{\theta}}{\partial\theta^{i}}\frac{\partial f_{\theta}}{\partial\theta^{j}}\Bigr],\label{eq:g_ij_derived_gen}
\end{equation}
where
\[
E_{\theta}''[\cdot]=\frac{\int_{T}(\cdot)\varphi''(f_{\theta})d\mu}{\int_{T}u_{0}\varphi'(f_{\theta})d\mu}.
\]
Because the operations of integration with respect to $\mu$ and differentiation
with respect to $\theta^{i}$ are commutative, we have
\begin{equation}
0=\frac{\partial}{\partial\theta^{i}}\int_{T}p_{\theta}d\mu=\int_{T}\frac{\partial}{\partial\theta^{i}}\varphi(f_{\theta})d\mu=\int_{T}\frac{\partial f_{\theta}}{\partial\theta^{i}}\varphi'(f_{\theta})d\mu,\label{eq:score_step_gen}
\end{equation}
and 
\begin{equation}
0=\int_{T}\frac{\partial^{2}f_{\theta}}{\partial\theta^{i}\partial\theta^{j}}\varphi'(f_{\theta})d\mu+\int_{T}\frac{\partial f_{\theta}}{\partial\theta^{i}}\frac{\partial f_{\theta}}{\partial\theta^{j}}\varphi''(f_{\theta})d\mu.\label{eq:second_derivative_gen}
\end{equation}
Thus expression \prettyref{eq:g_ij_derived_gen} follows from \prettyref{eq:second_derivative_gen}.
In addition, expression \prettyref{eq:score_step_gen} implies 
\begin{equation}
E_{\theta}'\Bigl[\frac{\partial f_{\theta}}{\partial\theta^{i}}\Bigr]=0.\label{eq:score_expectation_zero_gen}
\end{equation}

A consequence of \prettyref{eq:g_ij_derived_gen} is the correspondence
between the functions $\partial f_{\theta}/\partial\theta^{i}$ and
the basis vectors $\partial/\partial\theta^{i}$. The inner product
of vectors
\[
X=\sum_{i}a^{i}\frac{\partial}{\partial\theta^{i}}\qquad\text{and}\qquad Y=\sum_{i}b^{j}\frac{\partial}{\partial\theta^{j}}
\]
can be written as
\begin{equation}
g(X,Y)=\sum_{i,j}g_{ij}a^{i}b^{j}=\sum_{i,j}E_{\theta}''\Bigl[\frac{\partial f_{\theta}}{\partial\theta^{i}}\frac{\partial f_{\theta}}{\partial\theta^{j}}\Bigr]a^{i}b^{j}=E_{\theta}''[\widetilde{X}\widetilde{Y}],\label{eq:g_functional}
\end{equation}
where
\[
\widetilde{X}=\sum_{i}a^{i}\frac{\partial f_{\theta}}{\partial\theta^{i}}\qquad\text{and}\qquad\widetilde{Y}=\sum_{i}b^{j}\frac{\partial f_{\theta}}{\partial\theta^{j}}.
\]
As a result, the tangent space $T_{p_{\theta}}\mathcal{P}$ can be
identified with $\widetilde{T}_{p_{\theta}}\mathcal{P}$, which is
defined as the vector space spanned by $\partial f_{\theta}/\partial\theta^{i}$,
equipped with the inner product $\langle\widetilde{X},\widetilde{Y}\rangle_{\theta}=E_{\theta}''[\widetilde{X}\widetilde{Y}]$.
By \prettyref{eq:score_expectation_zero_gen}, if a vector $\widetilde{X}$
belongs to $\widetilde{T}_{p_{\theta}}\mathcal{P}$, then $E_{\theta}'[\widetilde{X}]=0$.
Independent of the definition of $(g_{ij})$, the expression in the
right-hand side of \prettyref{eq:g_functional} always defines a semi-inner
product in $\widetilde{T}_{p_{\theta}}\mathcal{P}$.

\subsection{Parametric $\varphi$-Families of Probability Distribution\label{sec:phi-family}}

Let $c\colon T\rightarrow\mathbb{R}$ be a measurable function such
that $p:=\varphi(c)$ is probability density in $\mathcal{P}_{\mu}$.
We take any measurable functions $u_{1},\dots u_{n}\colon T\rightarrow\mathbb{R}$
satisfying the following conditions:
\begin{itemize}
\item [(i)] $\int_{T}u_{i}\varphi'(c)d\mu=0$, and
\item [(ii)] there exists $\varepsilon>0$ such that
\[
\int_{T}\varphi(c+\lambda u_{i})d\mu<\infty,\qquad\text{for all }\lambda\in(-\varepsilon,\varepsilon).
\]

\end{itemize}
Define $\Theta\subseteq\mathbb{R}^{n}$ as the set of all vectors
$\theta=(\theta^{i})\in\mathbb{R}^{n}$ such that
\[
\int_{T}\varphi\biggl(c+\lambda\sum_{k=1}^{n}\theta^{i}u_{i}\biggr)d\mu<\infty,\qquad\text{for some }\lambda>1.
\]
The elements of the \textit{parametric $\varphi$-family} $\mathcal{F}_{p}=\{p(t;\theta):\theta\in\Theta\}$
centered at $p=\varphi(c)$ are given by the one-to-one mapping
\begin{equation}
p(t;\theta):=\varphi\biggl(c(t)+\sum_{i=1}^{n}\theta^{i}u_{i}(t)-\psi(\theta)u_{0}(t)\biggr),\qquad\text{for each }\theta=(\theta^{i})\in\Theta.\label{eq:phi_c}
\end{equation}
where the \textit{normalizing function} $\psi\colon\Theta\rightarrow[0,\infty)$
is introduced so that expression \prettyref{eq:phi_c} defines a probability
distribution in $\mathcal{P}_{\mu}$. 

Condition (ii) is always satisfied if the set $T$ is finite. It can
be shown that the normalizing function $\psi$ is also convex (and
the set $\Theta$ is open and convex). Under conditions (i)--(ii),
the family $\mathcal{F}_{p}$ is a submanifold of a non-parametric
$\varphi$-family. For the non-parametric case, we refer to \cite{Vigelis:2013,Vigelis:2013a}. 

By the equalities 
\[
\frac{\partial f_{\theta}}{\partial\theta^{i}}=u_{i}(t)-\frac{\partial\psi}{\partial\theta^{i}},\qquad-\frac{\partial^{2}f_{\theta}}{\partial\theta^{i}\partial\theta^{j}}=-\frac{\partial^{2}\psi}{\partial\theta^{i}\partial\theta^{j}},
\]
we get
\[
g_{ij}=\frac{\partial^{2}\psi}{\partial\theta^{i}\partial\theta^{j}}.
\]
In other words, the matrix $(g_{ij})$ is the Hessian of the normalizing
function $\psi$.

For $\varphi(\cdot)=\exp(\cdot)$ and $u_{0}=1$, expression \prettyref{eq:phi_c}
defines a parametric exponential family of probability\textit{ }distributions
\textit{$\mathcal{E}_{p}$}. In exponential families, the normalizing
function is recognized as the Kullback--Leibler divergence between
$p(t)$ and $p(t;\theta)$. Using this result, we can define the $\varphi$-divergence
$\mathcal{D}_{\varphi}(\cdot\mathbin{\Vert}\cdot)$, which generalizes
the Kullback--Leibler divergence $\mathcal{D}_{\mathrm{KL}}(\cdot\mathbin{\Vert}\cdot)$. 

By \prettyref{eq:phi_c} we can write 
\[
\psi(\theta)u_{0}(t)=\sum_{i=1}^{n}\theta^{i}u_{i}(t)+\varphi^{-1}(p(t))-\varphi^{-1}(p(t;\theta)).
\]
From condition (i), this equation yields 
\[
\psi(\theta)\int_{T}u_{0}\varphi'(c)d\mu=\int_{T}[\varphi^{-1}(p)-\varphi^{-1}(p_{\theta})]\varphi'(c)d\mu.
\]
In view of $\varphi'(c)=1/(\varphi^{-1})'(p)$, we get
\begin{equation}
\psi(\theta)=\frac{{\displaystyle \int_{T}\frac{\varphi^{-1}(p)-\varphi^{-1}(p_{\theta})}{(\varphi^{-1})'(p)}d\mu}}{{\displaystyle \int_{T}\frac{u_{0}}{(\varphi^{-1})'(p)}d\mu}}=:\mathcal{D}_{\varphi}(p\mathbin{\Vert}p_{\theta}),\label{eq:D_phi}
\end{equation}
which defines the $\varphi$-divergence $\mathcal{D}_{\varphi}(p\mathbin{\Vert}p_{\theta})$.
Clearly, expression \prettyref{eq:D_phi} can be used to extend the
definition of $\mathcal{D}_{\varphi}(\cdot\mathbin{\Vert}\cdot)$
to any probability distributions $p$ and $q$ in $\mathcal{P}_{\mu}$.

\section{$\alpha$-Connections\label{sec:alpha_connections}}

We use the $\varphi$-divergence $\mathcal{D}_{\varphi}(\cdot\mathbin{\Vert}\cdot)$
to define a pair of mutually dual connection in generalized statistical
manifolds. Let $\mathcal{D}\colon M\times M\rightarrow[0,\infty)$
be a non-negative, differentiable function defined on a smooth manifold
$M$, such that
\begin{equation}
\mathcal{D}(p\mathbin{\Vert}q)=0\quad\text{if and only if }\quad p=q.\label{eq:divergence_zero_equality}
\end{equation}
The function $\mathcal{D}(\cdot\mathbin{\Vert}\cdot)$ is called a
\emph{divergence} if the matrix $(g_{ij})$, whose entries are given
by
\begin{equation}
g_{ij}(p)=-\Bigl[\Bigl(\frac{\partial}{\partial\theta^{i}}\Bigr)_{p}\Bigl(\frac{\partial}{\partial\theta^{j}}\Bigr)_{q}\mathcal{D}(p\mathbin{\Vert}q)\Bigr]_{q=p},\label{eq:gij_divergence_first}
\end{equation}
is positive definite for each $p\in M$. Hence a divergence $\mathcal{D}(\cdot\mathbin{\Vert}\cdot)$
defines a metric in $M$. A divergence $\mathcal{D}(\cdot\mathbin{\Vert}\cdot)$
also induces a pair of mutually dual connections $D$ and $D^{*}$,
whose Christoffel symbols are given by
\begin{equation}
\Gamma_{ijk}=-\Bigl[\Bigl(\frac{\partial^{2}}{\partial\theta^{i}\partial\theta^{j}}\Bigr)_{p}\Bigl(\frac{\partial}{\partial\theta^{k}}\Bigr)_{q}\mathcal{D}(p\mathbin{\Vert}q)\Bigr]_{q=p}\label{eq:Gamma_ijk_divergence}
\end{equation}
and 
\begin{equation}
\Gamma_{ijk}^{*}=-\Bigl[\Bigl(\frac{\partial}{\partial\theta^{k}}\Bigr)_{p}\Bigl(\frac{\partial^{2}}{\partial\theta^{i}\partial\theta^{j}}\Bigr)_{q}\mathcal{D}(p\mathbin{\Vert}q)\Bigr]_{q=p},\label{eq:Gamma_star_ijk_divergence}
\end{equation}
respectively. By a simple computation, we get 
\[
\frac{\partial g_{jk}}{\partial\theta^{i}}=\Gamma_{ijk}+\Gamma_{ikj}^{*},
\]
showing that $D$ and $D^{*}$ are mutually dual.

In \prettyref{sec:phi-family}, the $\varphi$-divergence between
two probability distributions $p$ and $q$ in $\mathcal{P}_{\mu}$
was defined as
\begin{equation}
\mathcal{D}_{\varphi}(p\mathbin{\Vert}q):=\frac{{\displaystyle \int_{T}\frac{\varphi^{-1}(p)-\varphi^{-1}(q)}{(\varphi^{-1})'(p)}d\mu}}{{\displaystyle \int_{T}\frac{u_{0}}{(\varphi^{-1})'(p)}d\mu}}.\label{eq:P_phi_pq}
\end{equation}
Because $\varphi$ is convex, it follows that $\mathcal{D}_{\varphi}(p\mathbin{\Vert}q)\geq0$
for all $p,q\in\mathcal{P}_{\mu}$. In addition, if we assume that
$\varphi(\cdot)$ is strictly convex, then $\mathcal{D}_{\varphi}(p\mathbin{\Vert}q)=0$
if and only if $p=q$. In a generalized statistical manifold $\mathcal{P}=\{p(t;\theta):\theta\in\Theta\}$,
the metric derived from the divergence $\mathcal{D}(q\mathbin{\Vert}p):=\mathcal{D}_{\varphi}(p\mathbin{\Vert}q)$
coincides with \prettyref{eq:g_ij_gen}. Expressing the $\varphi$-divergence
$\mathcal{D}_{\varphi}(\cdot\mathbin{\Vert}\cdot)$ between $p_{\theta}$
and $p_{\vartheta}$ as 
\[
\mathcal{D}(p_{\theta}\mathbin{\Vert}p_{\vartheta})=E_{\vartheta}'[(f_{\vartheta}-f_{\theta})],
\]
after some manipulation, we get
\begin{align*}
g_{ij} & =-\Bigl[\Bigl(\frac{\partial}{\partial\theta^{i}}\Bigr)_{p}\Bigl(\frac{\partial}{\partial\theta^{j}}\Bigr)_{q}\mathcal{D}(p\mathbin{\Vert}q)\Bigr]_{q=p}\\
 & =-E_{\theta}'\Bigl[\frac{\partial^{2}f_{\theta}}{\partial\theta^{i}\partial\theta^{j}}\Bigr].
\end{align*}
As a consequence, expression \prettyref{eq:P_phi_pq} defines a divergence
on statistical manifolds. 

Let $D^{(1)}$ and $D^{(-1)}$ denote the pair of dual connections
derived from $\mathcal{D}_{\varphi}(\cdot\mathbin{\Vert}\cdot)$.
By \prettyref{eq:Gamma_ijk_divergence} and \prettyref{eq:Gamma_star_ijk_divergence},
the Christoffel symbols $\Gamma_{ijk}^{(1)}$ and $\Gamma_{ijk}^{(-1)}$
are given by 
\begin{equation}
\Gamma_{ijk}^{(1)}=E_{\theta}''\Bigl[\frac{\partial^{2}f_{\theta}}{\partial\theta^{i}\partial\theta^{j}}\frac{\partial f_{\theta}}{\partial\theta^{k}}\Bigr]-E_{\theta}'\Bigl[\frac{\partial^{2}f_{\theta}}{\partial\theta^{i}\partial\theta^{j}}\Bigr]E_{\theta}''\Bigl[u_{0}\frac{\partial f_{\theta}}{\partial\theta^{k}}\Bigr]\label{eq:Gamma_one}
\end{equation}
and 
\begin{align}
\Gamma_{ijk}^{(-1)} & =E_{\theta}''\Bigl[\frac{\partial^{2}f_{\theta}}{\partial\theta^{i}\partial\theta^{j}}\frac{\partial f_{\theta}}{\partial\theta^{k}}\Bigr]+E_{\theta}'''\Bigl[\frac{\partial f_{\theta}}{\partial\theta^{i}}\frac{\partial f_{\theta}}{\partial\theta^{j}}\frac{\partial f_{\theta}}{\partial\theta^{k}}\Bigr]\nonumber \\
 & \qquad-E_{\theta}''\Bigl[\frac{\partial f_{\theta}}{\partial\theta^{j}}\frac{\partial f_{\theta}}{\partial\theta^{k}}\Bigr]E_{\theta}''\Bigl[u_{0}\frac{\partial f_{\theta}}{\partial\theta^{i}}\Bigr]-E_{\theta}''\Bigl[\frac{\partial f_{\theta}}{\partial\theta^{i}}\frac{\partial f_{\theta}}{\partial\theta^{k}}\Bigr]E_{\theta}''\Bigl[u_{0}\frac{\partial f_{\theta}}{\partial\theta^{j}}\Bigr],\label{eq:Gamma_minus_one}
\end{align}
where
\[
E_{\theta}'''[\cdot]=\frac{\int_{T}(\cdot)\varphi'''(f_{\theta})d\mu}{\int_{T}u_{0}\varphi'(f_{\theta})d\mu}.
\]
Notice that in parametric $\varphi$-families, the Christoffel symbols
$\Gamma_{ijk}^{(1)}$ vanish identically. Thus, in these families,
the connection $D^{(1)}$ is flat.

Using the pair of mutually dual connections $D^{(1)}$ and $D^{(-1)}$,
we can specify a family of $\alpha$-connections $D^{(\alpha)}$ in
generalized statistical manifolds. The Christoffel symbol of $D^{(\alpha)}$
is defined by
\begin{equation}
\Gamma_{ijk}^{(\alpha)}=\frac{1+\alpha}{2}\Gamma_{ijk}^{(1)}+\frac{1-\alpha}{2}\Gamma_{ijk}^{(-1)}.\label{eq:alpha-connection}
\end{equation}
The connections $D^{(\alpha)}$ and $D^{(-\alpha)}$ are mutually
dual, since
\[
\frac{\partial g_{jk}}{\partial\theta^{i}}=\Gamma_{ijk}^{(\alpha)}+\Gamma_{ikj}^{(-\alpha)}.
\]
For $\alpha=0$ , the connection $D^{(0)}$, which is clearly self-dual,
corresponds to the Levi--Civita connection $\nabla$. One can show
that $\Gamma_{ijk}^{(0)}$ can be derived from the expression defining
the Christoffel symbols of $\nabla$ in terms of the metric: 
\[
\Gamma_{ijk}=\sum_{m}\Gamma_{ij}^{m}g_{mk}=\frac{1}{2}\Bigl(\frac{\partial g_{ki}}{\partial\theta^{j}}+\frac{\partial g_{kj}}{\partial\theta^{i}}-\frac{\partial g_{ij}}{\partial\theta^{k}}\Bigr).
\]
The connection $D^{(\alpha)}$ can be equivalently defined by 
\[
\Gamma_{ijk}^{(\alpha)}=\Gamma_{ijk}^{(0)}-\alpha T_{ijk},
\]
where
\begin{multline}
T_{ijk}=\frac{1}{2}E_{\theta}'''\Bigl[\frac{\partial f_{\theta}}{\partial\theta^{i}}\frac{\partial f_{\theta}}{\partial\theta^{j}}\frac{\partial f_{\theta}}{\partial\theta^{k}}\Bigr]-\frac{1}{2}E_{\theta}''\Bigl[\frac{\partial f_{\theta}}{\partial\theta^{k}}\frac{\partial f_{\theta}}{\partial\theta^{i}}\Bigr]E_{\theta}''\Bigl[u_{0}\frac{\partial f_{\theta}}{\partial\theta^{j}}\Bigr]\\
-\frac{1}{2}E_{\theta}''\Bigl[\frac{\partial f_{\theta}}{\partial\theta^{k}}\frac{\partial f_{\theta}}{\partial\theta^{j}}\Bigr]E_{\theta}''\Bigl[u_{0}\frac{\partial f_{\theta}}{\partial\theta^{i}}\Bigr]-\frac{1}{2}E_{\theta}''\Bigl[\frac{\partial f_{\theta}}{\partial\theta^{i}}\frac{\partial f_{\theta}}{\partial\theta^{j}}\Bigr]E_{\theta}''\Bigl[u_{0}\frac{\partial f_{\theta}}{\partial\theta^{k}}\Bigr].\label{eq:Tijk}
\end{multline}
 In the case that $\varphi$ is the exponential function and $u_{0}=1$,
equations \prettyref{eq:Gamma_one}, \prettyref{eq:Gamma_minus_one},
\prettyref{eq:alpha-connection} and \prettyref{eq:Tijk} reduce to
the classical expressions for statistical manifolds.

\subsection{Parallel Transport\label{sec:parallel_transport}}

Let $\gamma\colon I\rightarrow M$ be a smooth curve in a smooth manifold
$M$, with a connection $D$. A vector field $V$ along $\gamma$
is said to be \emph{parallel} if $D_{d/dt}V(t)=0$ for all $t\in I$.
Take any tangent vector $V_{0}$ at $\gamma(t_{0})$, for some $t_{0}\in I$.
Then there exists a unique vector field $V$ along $\gamma$, called
the \emph{parallel transport} of $V_{0}$ along $\gamma$, such that
$V(t_{0})=V_{0}$. 

A connection $D$ can be recovered from the parallel transport. Fix
any smooth vectors fields $X$ and $Y$. Given $p\in M$, define $\gamma\colon I\rightarrow M$
to be an integral curve of $X$ passing through $p$. In other words,
$\gamma(t_{0})=p$ and $\frac{d\gamma}{dt}=X(\gamma(t))$. Let $P_{\gamma,t_{0},t}\colon T_{\gamma(t_{0})}M\rightarrow T_{\gamma(t)}M$
denote the parallel transport of a vector along $\gamma$ from $t_{0}$
to $t$. Then we have
\[
(D_{X}Y)(p)=\frac{d}{dt}P_{\gamma,t_{0},t}^{-1}(Y(c(t))\biggr|_{t=t_{0}}.
\]
For details, we refer to \cite{Carmo:2013}.

To avoid some technicalities, we assume that the set $T$ is finite.
In this case, we can consider a generalized statistical manifold $\mathcal{P}=\{p(t;\theta):\theta\in\Theta\}$
for which $\mathcal{P}=\mathcal{P}_{\mu}$. The connection $D^{(1)}$
can be derived from the parallel transport 
\[
P_{q,p}\colon\widetilde{T}_{q}\mathcal{P}\rightarrow\widetilde{T}_{p}\mathcal{P}
\]
given by
\[
\widetilde{X}\mapsto\widetilde{X}-E_{\theta}'[\widetilde{X}]u_{0},
\]
where $p=p_{\theta}$. Recall that the tangent space $T_{p}\mathcal{P}$
can be identified with $\widetilde{T}_{p}\mathcal{P}$, the vector
space spanned by the functions $\partial f_{\theta}/\partial\theta^{i}$,
equipped with the inner product $\langle\widetilde{X},\widetilde{Y}\rangle=E_{\theta}''[\widetilde{X}\widetilde{Y}]$,
where $p=p_{\theta}$. We remark that $P_{q,p}$ does not depend on
the curve joining $q$ and $p$. As a result, the connection $D^{(1)}$
is flat. Denote by $\gamma(t)$ the coordinate curve given locally
by $\theta(t)=(\theta^{1},\dots,\theta^{i}+t,\dots,\theta^{n})$.
Observing that $P_{\gamma(0),\gamma(t)}^{-1}$ maps the vector $\frac{\partial f_{\theta}}{\partial\theta^{j}}(t)$
to 
\[
\frac{\partial f_{\theta}}{\partial\theta^{j}}(t)-E_{\theta(0)}'\Bigl[\frac{\partial f_{\theta}}{\partial\theta^{j}}(t)\Bigr]u_{0},
\]
we define the connection
\begin{align*}
\widetilde{D}_{\partial f_{\theta}/\partial\theta_{i}}\frac{\partial f_{\theta}}{\partial\theta_{j}} & =\frac{d}{dt}P_{\gamma(0),\gamma(t)}^{-1}\Bigl(\frac{\partial f_{\theta}}{\partial\theta_{j}}(\gamma(t)\Bigr)\biggr|_{t=0}\\
 & =\frac{d}{dt}\Bigl(\frac{\partial f_{\theta(t)}}{\partial\theta^{j}}-E_{\theta(0)}'\Bigl[\frac{\partial f_{\theta(t)}}{\partial\theta^{j}}\Bigr]u_{0}\Bigr)\biggr|_{t=0}\\
 & =\frac{\partial^{2}f_{\theta}}{\partial\theta^{i}\partial\theta^{j}}-E_{\theta}'\Bigl[\frac{\partial^{2}f_{\theta}}{\partial\theta^{i}\partial\theta^{j}}\Bigr]u_{0}.
\end{align*}
Let us denote by $D$ the connection corresponding to $\widetilde{D}$,
which acts on smooth vector fields in $T_{p}\mathcal{P}$. By this
identification, we have
\begin{align*}
g\Bigl(D_{\partial/\partial\theta_{i}}\frac{\partial}{\partial\theta_{j}},\frac{\partial}{\partial\theta_{k}}\Bigr) & =\Bigl\langle\widetilde{D}_{\partial f_{\theta}/\partial\theta_{i}}\frac{\partial f_{\theta}}{\partial\theta_{j}},\frac{\partial f_{\theta}}{\partial\theta_{k}}\Bigr\rangle\\
 & =E_{\theta}''\Bigl[\frac{\partial^{2}f_{\theta}}{\partial\theta^{i}\partial\theta^{j}}\frac{\partial f_{\theta}}{\partial\theta^{k}}\Bigr]-E_{\theta}'\Bigl[\frac{\partial^{2}f_{\theta}}{\partial\theta^{i}\partial\theta^{j}}\Bigr]E_{\theta}''\Bigl[u_{0}\frac{\partial f_{\theta}}{\partial\theta^{k}}\Bigr]\\
 & =\Gamma_{ijk}^{(1)},
\end{align*}
showing that $D=D^{(1)}$.

\bibliographystyle{plain}
\bibliography{refs}

\def\cprime{$'$}
\begin{thebibliography}{10}

\bibitem{Amari:2000}
Shun-ichi Amari and Hiroshi Nagaoka.
\newblock {\em Methods of information geometry}, volume 191 of {\em
  Translations of Mathematical Monographs}.
\newblock American Mathematical Society, Providence, RI; Oxford University
  Press, Oxford, 2000.
\newblock Translated from the 1993 Japanese original by Daishi Harada.

\bibitem{Amari:2011}
Shun-ichi Amari and Atsumi Ohara.
\newblock Geometry of {$q$}-exponential family of probability distributions.
\newblock {\em Entropy}, 13(6):1170--1185, 2011.

\bibitem{Amari:2012}
Shun-ichi Amari, Atsumi Ohara, and Hiroshi Matsuzoe.
\newblock Geometry of deformed exponential families: invariant, dually-flat and
  conformal geometries.
\newblock {\em Phys. A}, 391(18):4308--4319, 2012.

\bibitem{Calin:2014}
Ovidiu Calin and Constantin Udri\c{s}te.
\newblock {\em Geometric Modeling in Probability and Statistics}.
\newblock Springer, 2014.

\bibitem{Carmo:2013}
Mafredo~P. do~Carmo.
\newblock {\em Riemannian Geometry}.
\newblock Birkh\"{a}user, 14th edition, 2013.

\bibitem{Kaniadakis:2002}
G.~Kaniadakis.
\newblock Statistical mechanics in the context of special relativity.
\newblock {\em Phys. Rev. E (3)}, 66(5):056125, 17, 2002.

\bibitem{Matsuzoe:2014}
Hiroshi Matsuzoe.
\newblock Hessian structures on deformed exponential families and their
  conformal structures.
\newblock {\em Differential Geom. Appl.}, 35(suppl.):323--333, 2014.

\bibitem{Shima:2007}
Hirohiko Shima.
\newblock {\em The geometry of {H}essian structures}.
\newblock World Scientific Publishing Co. Pte. Ltd., Hackensack, NJ, 2007.

\bibitem{Tsallis:1994}
Constantino Tsallis.
\newblock What are the numbers that experiments provide?
\newblock {\em Quimica Nova}, 17(6):468--471, 1994.

\bibitem{Vigelis:2013a}
Rui~F. Vigelis and Charles~C. Cavalcante.
\newblock The {$\Delta_2$}-condition and {$\varphi$}-families of probability
  distributions.
\newblock In {\em Geometric science of information}, volume 8085 of {\em
  Lecture Notes in Comput. Sci.}, pages 729--736. Springer, Heidelberg, 2013.

\bibitem{Vigelis:2013}
Rui~F. Vigelis and Charles~C. Cavalcante.
\newblock On {$\varphi$}-families of probability distributions.
\newblock {\em J. Theoret. Probab.}, 26(3):870--884, 2013.

\end{thebibliography}

\end{document}